\documentclass[12pt,fleqn, a4]{article}

\usepackage{amsfonts}
\usepackage{stmaryrd}
\usepackage{amssymb}
\usepackage{amsthm}

\usepackage{amsfonts}
\usepackage{stmaryrd}
\usepackage{amssymb}
\usepackage{latexsym}
\usepackage{graphics}
\usepackage{amsthm}
\usepackage{latexsym,amsmath,amscd}
\usepackage{cleveref}
\usepackage{color}
\usepackage{dsfont}

\newcommand{\colval}{0.3}
\definecolor{colone}{gray}{\colval}

\newcommand{\dcb}{\begin{array}{lll}}
\newcommand{\dce}{\end{array}}
\newcommand{\ebe}{\begin{enumerate}\setlength{\baselineskip}{13pt}\setlength{\parskip}{5pt}}
\newcommand{\dbe}{\end{enumerate}}

\newcommand{\ibegin}{\begin{itemize}\setlength{\baselineskip}{19pt}\setlength{\parskip}{7pt}}
\newcommand{\iend}{\end{itemize}}
\newcommand{\ok}{\rule{4pt}{6pt}}
\newcommand{\desb}{\begin{description}}
\newcommand{\dese}{\end{description}}

\newtheorem{Theorem}{Theorem}[section]
\newtheorem {Cor}[Theorem]{Corollary}
\newtheorem {definition}[Theorem]{Definition}
\newtheorem {pro}[Theorem]{Proposition}
\newtheorem {Lemma}[Theorem]{Lemma}
\newtheorem {rem}[Theorem]{Remark}
\newtheorem {assumption}[Theorem]{Assumption}

\newcommand {\bd}{\begin{definition}}
\newcommand {\ed}{\end{definition}}
\newcommand {\bpro}{\begin{pro}}
\newcommand {\epro}{\end{pro}}
\newcommand {\bl}{\begin{Lemma}}
\newcommand {\el}{\end{Lemma}}
\newcommand {\bcor}{\begin{Cor}}
\newcommand {\ecor}{\end{Cor}}
\newcommand {\brem }{\begin{rem} \rm }
\newcommand {\erem }{\end{rem}}
\newcommand{\bethe}{\begin{Theorem}}
\newcommand{\ethe}{\end{Theorem}}
\newcommand {\bassumption}{\begin{assumption}}
\newcommand {\eassumption}{\end{assumption}}

\def \ind{1\!\!1}

\begin{document}

\begin{center}
\textbf{\Large Dynamic construction of martingales of density functions}\footnote{A preliminary version}
\end{center}

\begin{center}
Shiqi Song

{\footnotesize Laboratoire Analyse et Probabilités\\
Université d'Evry Val D'Essonne, France\\
shiqi.song@univ-evry.fr}
\end{center}

\begin{quote}
The density hypothesis on random times becomes now a standard in modeling of risks. One of the basic reasons to introduce the density hypothesis is the desire to have a computable credit risk model. However, recent work shows that merely an existence of a density function for the conditional law of the random times will not be enough for the purposes of some numerical implantation problems. It becomes necessary to have models with martingales of density functions evolving along with the development of the information flow, in particular, to have Markovian martingales of density functions determined by a stochastic differential equation. The quetion of constructing a martingale of density functions by a stochastic differential equation has been answered in one dimensional case. The aim of this note is to provide a solution in higher dimensional cases.  
\end{quote}

\section{Introduction}

Since the paper \cite{ejj,Pham} the density hypothesis on random times becomes a standard in the models involving risk. One of the basic reasons to introduce the density hypothesis is the desire to have a computable credit risk model (cf. \cite{CJZ,ejj, KL, Pham} for some computations and applications). However, the recent work \cite{wu} shows clearly that, in the absence of a simple dynamic equation for the martingale of the density functions (issued from the density hypothesis) with respect to the underlying filtration, the density functions can be inapplicable in some numerical implantation problem. More generally, for diverse reasons, one may hope to have the density functions evolving along with the development of the information flow, in particular, to have a Markovian martingale of density functions. This observaton raises the questions how to construct a martingale of density functions in a dynamic way, how to make this martingale evoluing in a predesigned way and how about the way the numerical simulation of such a martingale can be implanted.   

Note simply that any positive random variable is an stochastic exponential. Conversely, if $H_t(\omega,x)$ is a bounded predictable function in $x$ and if $Y$ is a martingale with bounded bracket $[Y,Y]$ and with the property $H_t(x)\Delta_t Y>-1$, the process $\mathcal{E}(H(x)\centerdot Y)_t, t\in\mathbb{R}_+, x\in[0,\infty]$, defines a non negative martingale such that, for any probability measure $\mu$ on $[0,\infty]$,$$
\mathbb{E}[\int_{[0,\infty]}\mathcal{E}(H(x)\centerdot Y)_t \mu(dx)]
=
\int_{[0,\infty]}\mathbb{E}[\mathcal{E}(H(x)\centerdot Y)_t] \mu(dx) =1.
$$ 
However, this is not enough, because, to be a density function for the conditional law of a random time, the process needs to satisfy a more restrictive equation
\begin{equation}\label{et=1}
\int_{[0,\infty]}\mathcal{E}(H(x)\centerdot Y)_t \mu(dx)=1
\end{equation}
for almost all $\omega$ and for all $t\in\mathbb{R}_+$. 
One may hope to modify dynamically the predictable function $H_t(x)$ to satisfy equation (\ref{et=1}). But, it happens that a better way to do so is to define dynamically a process $(C_t)_{t\in\mathbf{R}_+}$ of probability measures which are absolutely continuous with respect to $\mu$. That is the approach adopted in this note.

In one dimensional case, the question of dynamic construction of a martingale of density functions has an answer in this way by the papers \cite{JS2,songdynamic}. In these papers, a system of stochastic differential equations are introduced whose solutions define the distribution function of the conditional laws of a random time with respect to a given filtration. Differentiating the distribution functions will give the density functions. The derivatives of this kind of random distribution functions are computed in \cite{songdynamic} and a sufficient condition to have density functions in this way is that the Az\'ema supermartingale of the random time has an absolutely continuous drift (which is realized when the random time has an intensity).

The papers \cite{JS2,songdynamic} concerns only the case of a single random time. The question whether the program designed there can also be realized for a family of multiple random times, seems to have a negative answer. It is difficult to reproduce the methodology of \cite{JS2,songdynamic} in multiple random times cases. In this note we will follow a different approach : instead of constructing the distribution functions, we will directly construct probability valued processes via stochastic differential equation. This idea works with affirmative results. Along with this construction of probability valued processes we will exhibit the basic reasons for a solution of a stochastic differential equation to define probability measures. We will also discover the technical difficulties to prove the existence and uniqueness of the solution of such kind of stochastic differential equation. The discussion of this note gives some new insight on the problem of dynamical construction of martingales of density functions.

\

\section{Setup and preliminary analysis}

Let us fix the framework within which our study goes through. Let $\mathtt{X}$ be a compact metric space (for example, $\mathtt{X}=[0,T]^\mathfrak{l}$ for some real $T>0$ and some positive integer $\mathfrak{l}$). Let $\mathcal{M}(\mathtt{X})$ be the set of bounded signed Borel measures on $\mathtt{X}$. More generally, if $\mathcal{D}$ denotes a Borel sub-$\sigma$-algebra on $\mathtt{X}$, we denote by $\mathcal{M}(\mathtt{X},\mathcal{D})$ the set of all bounded signed measure on $\mathcal{D}$. 

We will consider processes taking values in $\mathcal{M}(\mathtt{X},\mathcal{D})$, especially, stochastic integrals of such processes. We adopt the following definition. Let $K_t(\omega,c,dx)$ be such that
\ebe
\item
For fixed $\omega\in\Omega$ and $t\geq 0$ and predictable process $\mathbf{c}$ taking values in $\mathcal{M}(\mathtt{X},\mathcal{D})$, $K_t(\omega,\mathbf{c},dx)$ is a bounded (signed) Borel measure on $\mathcal{D}$ taking values in $\mathbb{R}^d$,
\item
For any bounded $\mathcal{D}$ measurable function $h(x)$, for any predictable process $\mathbf{c}$ taking values in $\mathcal{M}(\mathtt{X},\mathcal{D})$, $\int_{\mathtt{X}} h(x)K_t(\cdot,\mathbf{c},dx)$ is a $d$-dimensional locally bounded predictable process.
\dbe
With this notion we define the measure valued stochastic differential equation. Let $Y$ be a $d$-dimensional local martingale and $\mathbf{c}$ a predictable $\mathcal{M}(\mathtt{X},\mathcal{D})$ valued process. We denote by $\int_0^t K_s(\mathbf{c})^\top dY_s, t\in\mathbb{R}_+$, the application from the set of all bounded $\mathcal{D}$ measurable functions into the space of local martingales: $$
h\longrightarrow \int_0^t \int_{\mathtt{X}}h(x)K_s(\mathbf{c}, dx)^\top dY_s,\ t\in\mathbb{R}_+.
$$ 
where the superscript $^{\top}$ designs the transposition of a vector. For a $\mathcal{F}_0$ measurable random measure $\mu$, for a process $C$ taking values in $\mathcal{M}(\mathtt{X},\mathcal{D})$, we consider the following stochastic differential equation:$$
\int_{\mathtt{X}}h(x)C_t(dx)
=\int_{\mathtt{X}}h(x)\mu(dx)
+\int_0^t \int_{\mathtt{X}}h(x)K_s(C_-,dx)^\top dY_s,
$$
for all bounded $\mathcal{D}$ measurable function $h$ on $\mathtt{X}$ and $t\in\mathbb{R}_+$. We use the following abbreviation to denote the above stochastic differential equation.
\begin{equation}\label{cequak}
\left\{
\dcb
dC_t = K_t(C_-)^\top dY_t,\\
C_0=\mu,
\dce
\right.
\end{equation}
The first question one should ask is naturally whether such an equation can have solutions. But what is equally important here is under what conditions the solutions of such an equation are probability valued. This amounts to check two conditions. When a solution gives rise a positive measure ? When the total mass of a solution remains to be 1 for all $t\in\mathbb{R}_+$. 

We adopt the notation$$
C_t[h]=\int_{\mathtt{X}}h(x)C_t(dx)
$$
for bounded $\mathcal{D}$ measurable function $h$. We make below some heuristic arguments. To have the positivity of a solution $C$ (when $\mu$ is positive), we can simply write$$
C_t[h]
=\int_{\mathtt{X}}h(x)\mu(dx)
+\int_0^t C_{s-}[h]\frac{\int_{\mathtt{X}}h(x)K_s(C_-,dx)^\top}{C_{s-}[h]}dY_s.
$$
Hence, to have the positivity, it is enough to assume the \textbf{condition I} :$$
\frac{\int_{\mathtt{X}}h(x)K_s(C_-,dx)}{C_{s-}[h]}
$$
is integrable with respect to $Y$ and $$
\frac{\int_{\mathtt{X}}h(x)K_s(C_-,dx)^\top}{C_{s-}[h]}\Delta_s Y>-1.
$$
Consider then the probability property (when $\mu$ is a probability measure) of a solution $C$, i.e. $C_t[1]\equiv 1$. Writting this condition with the stochastic differential equation, we get$$
1=1+\int_0^t \int_{\mathtt{X}}K_s(C_-,dx)^\top dY_s.
$$
A suitable condition would be \textbf{condition II} :$\int_{\mathtt{X}}K_s(C_-,dx)=0$. 

In this paper we will show that no trivial kernals $K$ satisfying the \textbf{condition I} and \textbf{condition II}  exist whose associated equations (\ref{cequak}) have probability valued solutions.

\

\section{Solutions of equation (\ref{cequak}) for finitely defined coefficient}\label{firstclass}

To carry out our programm we will concentrate ourself on the case where the coefficient $K$ in equation (\ref{cequak}) takes the form
\begin{equation}\label{kC}
K_s(\mathbf{c}, dx)=k_s(\mathbf{c},x)\mathbf{c}_{s}(dx),
\end{equation}
where $k_s(\omega,\mathbf{c},x)$ is a $\mathcal{P}(\mathbb{F}\otimes\mathcal{B}(\mathcal{M}(\mathtt{X},\mathcal{D})))\otimes\mathcal{D}$ measurable function bounded by a constant $\epsilon''$. 

In this section we begin the study of equation (\ref{cequak}) with a simple situation. We suppose that the $\sigma$-algebra $\mathcal{D}$ on $\mathtt{X}$ is generated by a finite partition $\mathfrak{P}$.

\subsection{Positivity of a solution}

Consider the \textbf{condition I} and the positivity of the solution of the equation (\ref{cequak}). 

\bl\label{cleqmu1}
Suppose that the function $k_s(\omega, \mathbf{c},x)$ is such that 
\begin{equation}\label{kDY>-1}
k_t(\mathbf{c},x)^\top\Delta_tY>-1
\end{equation} 
for all $\mathbf{c}$, $t$ and $x$. Let $C$ be a solution to equation (\ref{cequak}) on $\mathcal{D}$ with a positive initial measure $\mu$. Then, for any $A\in\mathfrak{P}$, $$
C_t(A)=\mu(A)\mathcal{E}(k(C_-,A)^\top\centerdot Y)_t,\ t\in\mathbb{R}_+,
$$
where $k_s(\mathbf{c},A)$ is the unique value of $k_s(\mathbf{c},x)$ for $x\in A$. In particular, $C_t(A)\geq 0$.
\el

\textbf{Proof.}
$$
\dcb
&&C_t(A)\\

&=&\mu(A)+\int_0^t \left(\int_{\mathtt{X}}\ind_A(x)k_s(C_-,x)^\top C_{s-}(dx)\right)dY_s\\

&=&\mu(A)+\int_0^t \left(k_s(C_-,A)^\top\int_{\mathtt{X}}\ind_A(x)C_{s-}(dx)\right)dY_s\\

&=&\mu(A)+\int_0^t k_s(C_-,A)^\top C_{s-}[A]dY_s.
\dce
$$
\ok

From now on in this paper we only consider kernal $K$ satisfying the condition (\ref{kC}) and (\ref{kDY>-1}).

\subsection{Probability property}

\bl\label{proba}
Let $C$ be the solution of the equation (\ref{cequak}) on $\mathcal{D}$. Suppose that $\mu$ is a probability measure and 
\begin{equation}\label{kernalC=0}
\int_{\mathtt{X}}k_s(C_-,x)^\top C_{s-}(dx)=0
\end{equation}
whenever $C_{s-}(\mathtt{X})\geq \frac{1}{3}$. Then, for any $t\in\mathbb{R}_+$, $C_t$ is a probability measure.
\el

\textbf{Proof.}
Let $R=\inf\{s\in\mathbb{R}_+: C_{s-}(\mathtt{X}) < \frac{1}{3}\}$. Clearly $R$ is strictly positive. We write$$
C_{t\wedge R}(\mathtt{X})
=\int_{\mathtt{X}}\mu(dx)
+\int_0^{t\wedge R}\int_{\mathtt{X}}k_s(C_-,x)^\top C_{s-}(dx)dY_s
=1.
$$
Notice that, if $R<t$, we have $1=C_{R}(\mathtt{X})\leq \frac{1}{3}$, a contradiction. \ok

\subsection{Existence of solutions on the finite $\sigma$-algebra $\mathcal{D}$}

For the existence problem of equation (\ref{cequak}) on a finite $\sigma$-algebra $\mathcal{D}$, we need to give more structure to the coefficient $k_s(\mathbf{c},x)$. Let $g_t(\omega,\mathbf{y},x)$ be a bounded $\mathcal{P}(\mathbb{F})\otimes\mathcal{B}(\mathbb{R}^k)\otimes\mathcal{D}$ measurable $d$-dimensional function. Let $\mathbf{h}=(h_i: 1\leq i\leq k)$ be a vector of (non random) Borel functions with $h_1\equiv 1$. For $c\in\mathcal{M}(\mathtt{X})$, let $$
c[\mathbf{h}]
=\left(\int_{\mathtt{X}}h_1(x)c(dx),\ldots\ldots,\int_{\mathtt{X}}h_k(x)c(dx)\right).
$$
We introduce the coefficient whose $e$-th component takes the form: 
\begin{equation}\label{kC2}
k_{e,s}(\omega, \mathbf{c},x)
=g_{e,s}(\omega,\mathbf{c}_s[\mathbf{h}],x)-\frac{\mathbf{c}_s[g_{e,s}(\omega,\mathbf{c}_s[\mathbf{h}])]}{\frac{1}{3}\vee \mathbf{c}_s(\mathtt{X})},
\end{equation}
for $s\in\mathbb{R}_+, \omega\in\Omega, x\in\mathtt{X}$, and $\mathbf{c}$ predictable process taking values in $\mathcal{M}(\mathtt{X})$, where$$
\mathbf{c}_s[g_{e,s}(\omega,\mathbf{c}_s[\mathbf{h}])]
=
\int_{\mathtt{X}}g_{e,s}(\omega,\mathbf{c}_s[\mathbf{h}],x)\mathbf{c}_s(dx).
$$ 
The main reason of the particular definition (\ref{kC2}) is the following identity
$$
\int_{\mathtt{X}}g_{e,s}(\omega,c[\mathbf{h}],x)-\frac{c[g_{e,s}(\omega,c[\mathbf{h}])]}{\frac{1}{3}\vee c(\mathtt{X})}c(dx)=0
$$
satisfied by any $c\in\mathcal{M}(\mathtt{X})$ such that $c(\mathtt{X})\geq \frac{1}{3}$, which makes the coefficient $k_s(\mathbf{c},x)$ to satisfy condition (\ref{kernalC=0}). We suppose in the definition (\ref{kC2}) that the function $g$ is uniformly Lipschitzian in $\mathbf{y}$ and in $x$:$$
\dcb
|g_{e,t}(\omega,\mathbf{y},x)-g_{e,t}(\omega,\mathbf{y}',x)|
\leq \epsilon'\|\mathbf{y}-\mathbf{y}'\|=\epsilon'\sum_{i=1}^k|y_i-y'_i|,\\
|g_{e,t}(\omega,\mathbf{y},x)-g_{e,t}(\omega,\mathbf{y},x')|
\leq \epsilon'\|x-x'\|=\epsilon'\sum_{i=1}^d|x_i-x'_i|.
\dce
$$ 
for some constant $\epsilon'$. We also suppose that $g$ vanishes when $y_1$ runs outside of a fixed finite interval, say, $[0,5]$ (a condition to have a Lipschitian condition on the function $\tilde{k}$ in the lemma below). Clearly the above defined coefficient $k_s(\mathbf{c},x)$ is bounded. We denote its upper bound by $\epsilon''$. At last, we assume that $g$ is chosen to ensure condition (\ref{kDY>-1}):
$$
k_t(\omega,\mathbf{c},x)^\top\Delta_t Y(\omega)>-1.
$$
In sum we have now a coefficient $k_t(\mathbf{c},x)$ which is bounded Markovian Lipschtzian satisfying conditions (\ref{kDY>-1}) and (\ref{kernalC=0}). Below we also use the notation $k_s(\mathbf{h},\mathbf{c}, x)$ to design a coefficient defined by (\ref{kC2}).

\bl\label{finite-uniqueness}
Suppose that the initial measure $\mu$ is a probability measure. Suppose that the vector $\mathbf{h}$ in (\ref{kC2}) is composed of $\mathcal{D}$ measurable functions. With the coefficient satisfying condition (\ref{kC}) and (\ref{kC2}), equation (\ref{cequak}) has a unique solution taking probability values in $\mathcal{M}(\mathtt{X},\mathcal{D})$.
\el

\textbf{Proof.}
Denote $\mathbf{X}=(C(A), A\in\mathfrak{P})\in[0,1]^{k'}$, where the integer $k'$ denotes the number of elements in $\mathfrak{P}$. Then, there exists a family of $d$-dimensional $\mathcal{P}(\mathbb{F})\otimes\mathcal{B}[0,1]^{k'}$ measurable functions $\tilde{k}_s(\mathbf{x},B), B\in\mathfrak{P}$, such that the resolution of the equation (\ref{cequak}) is equivalent to the resolution of$$
\left\{
\dcb
dC_{t}(B) &= C_{t-}(B) \tilde{k}_t(\mathbf{X}_{t-},B)^\top  dY_t,\\
\\
C_{0}(B)&=\mu(B),
\dce
\right\}
\ \forall B\in\mathfrak{P}.
$$
Since the function $g$ is Lipschitzian in $\mathbf{y}\in\mathbb{R}^{k}$ and vanishes when $y_1$ is too big, the functions $\tilde{k}$ can be chosen Lipschitzian in $\mathbf{x}\in\mathbb{R}^{k'}$. According to \cite{protter}, the above equation has a unique solution. The probability property of the solution is the consequence of condition (\ref{kDY>-1}) and (\ref{kernalC=0}), satisfied by coefficient defined by (\ref{kC2}). \ok

\

\subsection{Extension of solutions on the whole Borel $\sigma$-algebra of $\mathtt{X}$}

We always assume the conditions (\ref{kC}), (\ref{kC2}) and (\ref{kDY>-1}) with a probability initial measure $\mu$.

\bl\label{projective}
Let $\mathcal{D}'$ be a second $\sigma$-algebra generated by a finite partition $\mathfrak{P}'$. Suppose $\mathfrak{P}\subset \mathfrak{P}'$. Then, the coefficient in (\ref{kC2}) with the vector $\mathbf{h}$ composed of $\mathcal{D}$ measurable functions is well defined with respect to $\mathcal{D}'$ and equation (\ref{cequak}) has a unique solution on $\mathcal{D}'$ whose restriction on $\mathcal{D}$ coincides with the solution on $\mathcal{D}$. 
\el

\textbf{Proof.} The idea is simple. Let us illustrate it in the case where $\mathfrak{P}=\{A_1,A_2\}$ and $\mathfrak{P}'=\{A_{1,1},A_{1,2},A_2\}$ such that $A_{1,1}\cup A_{1,2}=A_{1}$. With the notation in Lemma \ref{finite-uniqueness}, the function $k$ in (\ref{kC2}) considered with respect to $\mathcal{D}'$ takes the form $$
k_s(C_{s-}[\mathbf{h}],x)
=\tilde{k}_t((C_{t-}(A_{1,1})+C_{t-}(A_{1,2}),C_{t-}(A_{2})),B)
$$
for $x\in B\in\mathfrak{P}$. Let $c_{1,1,t}=C_{t}(A_{1,1}), c_{1,2,t}=C_{t}(A_{1,2}), c_{2,t}=C_{t}(A_{2})$. We write the equation 
$$
\left\{
\dcb
dc_{1,1,t} &= c_{1,1,t-} \tilde{k}_t((c_{1,1,t-}+c_{1,2,t-},c_{2,t-}),A_1)^\top  dY_t,\\
dc_{1,2,t} &= c_{1,2,t-} \tilde{k}_t((c_{1,1,t-}+c_{1,2,t-},c_{2,t-}),A_1)^\top  dY_t,\\
dc_{2,t} &= c_{2,t-} \tilde{k}_t((c_{1,1,t-}+c_{1,2,t-},c_{2,t-}),A_2)^\top  dY_t,\\
\\
c_{1,1,0}&=\mu(A_{1,1}),\\
c_{1,2,0}&=\mu(A_{1,2}),\\
c_{2,0}&=\mu(A_{2}).
\dce
\right.
$$
In particuler,$$
\dcb
d(c_{1,1,t}+c_{1,2,t}) &=& (c_{1,1,t-}+c_{1,2,t-}) \tilde{k}_t((c_{1,1,t-}+c_{1,2,t-},c_{2,t-}),A_1)^\top  dY_t,\\
c_{1,1,0}+c_{1,2,0}&=&\mu[A_{1,1}]+\mu[A_{1,2}]=\mu[A_{1}].
\dce
$$
We see that the measure valued process on $\mathcal{D}$ associated with $(c_{1,1,t}+c_{1,2,t})$ and $c_{2,t}$ satisfies equation (\ref{cequak}) on $\mathcal{D}$. We conclude the lemma by uniqueness of solution.
\ok

\bethe
There exists a unique process $C$ probability measure valued on $\mathcal{B}(\mathtt{X})$ such that its restriction onto any finite $\sigma$-algebra $\mathcal{D}'$ containing $\mathcal{D}$ is a solution of the equation (\ref{cequak}) on $\mathcal{D}'$.
\ethe

\textbf{Proof.}
Let $\mathcal{D}_n$ be an increasing sequence of $\sigma$-algebras such that $\mathcal{D}_0=\mathcal{D}$ and $\vee_{n\in\mathbb{N}}\mathcal{D}_n=\mathcal{B}(\mathtt{X})$. Let $C$ be the solution of the equation (\ref{cequak}) on $\mathcal{D}_n$. Thanks to Lemma \ref{projective}, we can look the $C$ as the same process. Now, by measure extension theorem, $C$ is extended to be a process taking probability values (cf. Lemma \ref{proba}) in the space $\mathcal{M}(\mathtt{X})$. \ok

\brem
Let $\mathbf{c}$ denote the solution of equation (\ref{cequak}) on $\mathcal{D}$. Then, define, for any Borel set $A$ in $\mathtt{X}$,$$
C_t(A)=\sum_{B\in\mathfrak{P}:\mu[B]>0}\mathcal{E}(k(\mathbf{c}_-,B)^\top\centerdot Y)_t\mu(A\cap B),\ t\in\mathbb{R}_+.
$$
Following Lemma \ref{cleqmu1} we can check that this process $C$ is the solution of equation (\ref{cequak}) on $\mathcal{B}(\mathtt{X})$.
\erem

\

\section{Equation (\ref{cequak}) with continuous coefficient}\label{concoe}

The previous Section \ref{firstclass} proves the existence of the solution of equation (\ref{cequak}) when the kernal $K$ is assumed to depend on the parameters $\mathbf{c}$ and $x$ in a discrete way in order to maintain the computation within a finite rang. We try now to extend the existence result to equation (\ref{cequak}) with coefficient of a "continuous" nature. Our approach is the approximation of "continuous" case by the "discrete" cases. Estimations are needed to control the variation of the solutions of equation (\ref{cequak}) when the coefficient changes. This is the key point. 

\subsection{Continuous coefficient with separated variables}\label{inequality}

We consider equation (\ref{cequak}) with a kernal satisfying condition (\ref{kC}) and (\ref{kC2}) starting from an initial probability measure $\mu$. We suppose in addition that the function $g_t(\omega,\mathbf{y},x)$ in (\ref{kC2}) takes the form$$
g_{e,t}(\omega,\mathbf{y},x)
=\overline{g}_{e,t}(\omega,\mathbf{y})\check{g}_e(v(x)),
$$
where $v$ is a $\mathcal{D}$ measurable map from $\mathtt{X}$ into itself, $\check{g}_e(x)$ are continuous functions on $\mathtt{X}$, $\overline{g}_{e,t}(\omega,\mathbf{y})$ are $\mathcal{P}(\mathbb{F})\otimes\mathcal{B}(\mathbb{R}^k)$ measurable functions. Note that the function $g_{e,t}$ has now the variable $x$ separated from the other variables. Lipschitzian conditions with constant $\epsilon'$ and boundedness by constant $\epsilon''$ on $\overline{g}$ are assumed to satisfy the conditions in (\ref{kC2}). Notice that with such a coefficient,
$$
k_{e,s}(\omega, \mathbf{c},x)
=\overline{g}_{e,s}(\mathbf{c}_s[\mathbf{h}])\left(\check{g}_e(v(x))-\mathbf{c}_s[\check{g}_e\circ v]\right),
$$
for probability valued process $\mathbf{c}$. The dependence of the coefficient $k_{e,s}(\omega, \mathbf{c},x)$ on the functions $\mathbf{h}$ and $v$ will be essential in the following estimations. To indicate this dependence, we will use the notation $k_{e,s}(\omega, \mathbf{h}, \mathbf{c},v,x)$ (or $k_{e,s}(\omega, \mathbf{h}, \mathbf{c},x)$ if $v(x)=x$) instead of $k_{e,s}(\omega, \mathbf{c},x)$. Naturally we assume that the above functions are chosen to maintain condition (\ref{kDY>-1}).

Let $(f_i)_{i\in\mathbb{N}}$ be a countable family of continuous functions on $\mathtt{X}$,  uniformly bounded by 1, closed by multiplication, which generates a $\mathbb{Q}$-vector space dense in the space of continuous functions on $\mathtt{X}$ (cf. \cite[Theorem 81.3]{rw}). We assume that the functions in the vector $\mathbf{h}$ and the functions $\check{g}_e$ are in this family.

\subsection{Cancellation of a jump in a stochastic differential equation}

Before begining the estimation, we need a general result.

\bl\label{raccordement}
Let $\tilde{Y}$ be a semimartingale and $F$ is a locally bounded predictable functional in the sense of \cite{protter}. Let $R$ be a strictly positive stoppging time. Let $X$ be a solution of $$
X_t=X_0+\int_0^t F_s(X_-) d\tilde{Y}_s, \ t\in\mathbb{R}_+.
$$
Then, $X^{R-}$ is a solution of $$
X_t=X_0+\int_0^t F_s(X_-) d\tilde{Y}^{R-}_s.
$$
Conversely, let $Z$ be a solution of the stochastic differential equation$$
Z_t=X_0+\int_0^t F_s(Z_-) d\tilde{Y}^{R-}_s, \ t\in\mathbb{R}_+.
$$
Set $$
X=Z^R+F_R(Z_-)\Delta_R\tilde{Y}\ind_{[R,\infty)}.
$$
Then, $X$ is a solution of $$
X_t=X_0+\int_0^t F_s(X_-) d\tilde{Y}^R_s, \ t\in\mathbb{R}_+
$$
\el

\textbf{Proof.}
For the first part of the lemma, we write
$$
\dcb
X^{R-}
&=&X^R-\Delta_RX\ind_{[R,\infty)}
=X_0+F(X_-)\centerdot \tilde{Y}^R-F_R(X_-)\Delta_R\tilde{Y}\ind_{[R,\infty)}\\
&=&X_0+F(X_-)\centerdot \tilde{Y}^R-F(X_-)\centerdot(\Delta_R\tilde{Y}\ind_{[R,\infty)})\\
&=&X_0+F(X_-)\centerdot \tilde{Y}^{R-}.
\dce
$$
For the second part,
$$
\dcb
X
&=&Z^R+F_R(Z_-)\Delta_R\tilde{Y}\ind_{[R,\infty)}\\
&=&X_0+F(Z_-)\centerdot \tilde{Y}^{R-}+F_R(Z_-)\Delta_R\tilde{Y}\ind_{[R,\infty)}\\
&=&X_0+F(X_-)\centerdot \tilde{Y}^{R-}+F_R(X_-)\Delta_R\tilde{Y}\ind_{[R,\infty)}\\
&=&X_0+F(X_-)\centerdot \tilde{Y}^{R-}+F(X_-)\centerdot(\Delta_R\tilde{Y}\ind_{[R,\infty)})\\
&=&X_0+F(X_-)\centerdot \tilde{Y}^R.
\dce
$$\ok

\subsection{A (semi-)norm on the probability valued processes}

Let $\mathcal{D}$ be a finite $\sigma$-algebra. Let $\mathbf{h}'=(h'_i:1\leq i\leq k)$ and $\mathbf{h}''=(h''_i:1\leq i\leq k)$ be two vectors of $\mathcal{D}$ measurable functions. Let $v'$ and $v''$ be two $\mathcal{D}$ measurable maps from $\mathtt{X}$ into itself. These two families of functions inserted into the function $g$ of the form of subsection \ref{inequality} define two coefficients $k_t(\mathbf{h}',\mathbf{c},v',x)$ and $k_t(\mathbf{h}'',\mathbf{c},v'',x)$, which define in their turn the corresponding $C'$ and $C''$ solutions of the equation (\ref{cequak}) with the same initial probability measure $\mu$. We recall the following notations. For any bounded measurable function (random or parametered) $f(x)$ we denote$$
C'_{s-}[f]=\int_{\mathtt{X}}f(x) C'_{s-}(dx),
$$
and similarly $\mu[f], C''_{s-}[f]$. We also denote$$
\dcb
k'_s=k_s(\mathbf{h}',C'_-,v',x),\ \mbox{ and similarly $k''_s$},\\
\check{g}'(x)=\check{g}(v'(x)), \ \mbox{ and similarly $\check{g}''(x)$},\\
\overline{g}'=\overline{g}(C'_{s-}[\mathbf{h}']), \ \mbox{ and similarly $\overline{g}''$}.
\dce
$$

\bl\label{ineq1}
For any stopping time $R$, set $\|\mathbf{c}'-\mathbf{c}''\|_{R-}=\sup_{j\in\mathbb{N}}\mathbb{E}[\sup_{u< R}|\mathbf{c}'_{u}[f_j]-\mathbf{c}''_{u}[f_j]|]$. Let $a'$ be a constant in \cite[Chapter V section 2 Theorem 2]{protter}. Suppose that the stopping time $R$ satisfies $\|Y^{R-}_e\|_{H^\infty}\leq a$, for some real $a>0$ and for any $1\leq e\leq d$, where the norm $\|\cdot\|_{H^\infty}$ is defined in \cite[Chapter V section 2]{protter}. Then,
$$
\dcb
(1-da'a(2k\epsilon'+3\epsilon''))\|C'-C''\|_{R-}

&\leq&4da'a\epsilon'\sum_{i=1}^k\|h_i-h'_i\|_\infty
+2da'a\epsilon'\sum_{i=1}^k\|h'_i-h''_i\|_\infty\\
&&
+2a'a\epsilon''\sum_{e=1}^d\|\check{g}'_{e}-\check{g}_{e}\|_\infty
+2a'a\epsilon''\sum_{e=1}^d\|\check{g}'_{e}-\check{g}''_{e}\|_\infty,
\dce
$$
where for a function $u$ on $\mathtt{X}$, $\|u\|_\infty=\sup_{x\in\mathtt{X}}|u(x)|$.
\el

\textbf{Proof.}
For any $f_j$, by Lemma \ref{raccordement}, we have equation (\ref{cequak}) for $C'^{R-}$:$$
C'^{R-}_t[f_j]=\mu[f_j]+\int_0^t C'_{s-}[f_jk'^\top_s]dY^{R-}_s,
$$
and similarly for $C''^{R-}$. Hence,$$
\dcb
&&\mathbb{E}[\sup_{u< R}|C'_{u}[f_j]-C''_{u}[f_j]|]\\

&=&\mathbb{E}[\sup_{u< R}|C'^{R-}_{u}[f_i]-C''^{R-}_{u}[f_i]|]\\

&=&\mathbb{E}[\sup_{u<R}\left|\int_0^u (C'_{s-}[f_jk'^\top_s]
- C''_{s-}[f_jk''^\top_s]) dY^{R-}_s\right|].
\dce
$$
By \cite[Chapter V section 2 Theorem 2 and 3]{protter},$$
\dcb
&&\mathbb{E}[\sup_{u<R}\left|\int_0^u (C'_{s-}[f_jk'^\top_s]
- C''_{s-}[f_jk''^\top_s]) dY^{R-}_s\right|]\\

&\leq&\sum_{e=1}^d\mathbb{E}[\sup_{u<R}\left|\int_0^u (C'_{s-}[f_jk'_{e,s}]
- C''_{s-}[f_jk''_{e,s}]) dY^{R-}_{e,s}\right|]\\

&\leq&\sum_{e=1}^da'\mathbb{E}[\sup_{u<R}\left|C'_{u-}[f_jk'_{e,u}]
- C''_{u-}[f_jk''_{e,u}]\right|]\|Y^{R-}_{e}\|_{H^\infty}\\

&\leq&\sum_{e=1}^da'a\mathbb{E}[\sup_{u<R}\left|\overline{g}'_uC'_{u-}[f_j\check{g}'_{e}]
- \overline{g}''_uC''_{u-}[f_j\check{g}''_{e}]\right|] \\
&&+\sum_{e=1}^da'a\mathbb{E}[\sup_{u<R}\left|\overline{g}'_uC'_{u-}[f_j]C'_{u-}[\check{g}'_{e}]
- \overline{g}''_uC''_{u-}[f_j]C''_{u-}[\check{g}''_{e}]\right|] \\

&\leq&\sum_{e=1}^da'a\mathbb{E}[\sup_{u<R}\left|\overline{g}'_uC'_{u-}[f_j\check{g}'_{e}]
- \overline{g}''_uC'_{u-}[f_j\check{g}'_{e}]\right|] \\

&&+\sum_{e=1}^da'a\mathbb{E}[\sup_{u<R}\left|\overline{g}''_uC'_{u-}[f_j\check{g}'_{e}]
- \overline{g}''_uC''_{u-}[f_j\check{g}'_{e}]\right|] \\

&&+\sum_{e=1}^da'a\mathbb{E}[\sup_{u<R}\left|\overline{g}''_uC''_{u-}[f_j\check{g}'_{e}]
- \overline{g}''_uC''_{u-}[f_j\check{g}''_{e}]\right|]\\

&&+\sum_{e=1}^da'a\mathbb{E}[\sup_{u<R}\left|\overline{g}'_uC'_{u-}[f_j]C'_{u-}[\check{g}'_{e}]
- \overline{g}''_uC'_{u-}[f_j]C'_{u-}[\check{g}'_{e}]\right|] \\

&&+\sum_{e=1}^da'a\mathbb{E}[\sup_{u<R}\left|\overline{g}''_uC'_{u-}[f_j]C'_{u-}[\check{g}'_{e}]
- \overline{g}''_uC''_{u-}[f_j]C'_{u-}[\check{g}'_{e}]\right|] \\

&&+\sum_{e=1}^da'a\mathbb{E}[\sup_{u<R}\left|\overline{g}''_uC''_{u-}[f_j]C'_{u-}[\check{g}'_{e}]
- \overline{g}''_uC''_{u-}[f_j]C''_{u-}[\check{g}'_{e}]\right|] \\

&&+\sum_{e=1}^da'a\mathbb{E}[\sup_{u<R}\left|\overline{g}''_uC''_{u-}[f_j]C''_{u-}[\check{g}'_{e}]
- \overline{g}''_uC''_{u-}[f_j]C''_{u-}[\check{g}''_{e}]\right|].

\dce
$$
We have$$
\dcb
&&\mathbb{E}[\sup_{u<R}\left|\overline{g}'_uC'_{u-}[f_j\check{g}'_{e}]
- \overline{g}''_uC'_{u-}[f_j\check{g}'_{e}]\right|]\\

&\leq&\mathbb{E}[\sup_{u<R}\left|\overline{g}'_u
- \overline{g}''_u\right|]\\

&\leq&\mathbb{E}[\sup_{u<R}\epsilon'\left\|C'_{u-}[\mathbf{h}']-C''_{u-}[\mathbf{h}'']\right\|]\\

&\leq&\sum_{i=1}^k\epsilon'\mathbb{E}[\sup_{u<R}\left|C'_{u-}[h'_i]-C''_{u-}[h''_i]\right|]\\

&\leq&\sum_{i=1}^k\epsilon'\mathbb{E}[\sup_{u<R}\left|C'_{u-}[h'_i]-C''_{u-}[h'_i]\right|]
+\sum_{i=1}^k\epsilon'\mathbb{E}[\sup_{u<R}\left|C''_{u-}[h'_i]-C''_{u-}[h''_i]\right|]\\

&\leq&\sum_{i=1}^k\epsilon'\mathbb{E}[\sup_{u<R}\left|C'_{u-}[h_i]-C''_{u-}[h_i]\right|]\\
&&\sum_{i=1}^k\epsilon'\mathbb{E}[\sup_{u<R}\left|C'_{u-}[h'_i]-C'_{u-}[h_i]\right|]
+\sum_{i=1}^k\epsilon'\mathbb{E}[\sup_{u<R}\left|C''_{u-}[h_i]-C''_{u-}[h'_i]\right|]\\
&&+\sum_{i=1}^k\epsilon'\mathbb{E}[\sup_{u<R}\left|C''_{u-}[h'_i]-C''_{u-}[h''_i]\right|]\\

&\leq&\sum_{i=1}^k\epsilon'\|C'-C''\|_{R-}
+2\sum_{i=1}^k\epsilon'\|h_i-h'_i\|_\infty
+\sum_{i=1}^k\epsilon'\|h'_i-h''_i\|_\infty,\\
\dce
$$
and
$$
\dcb
&&\mathbb{E}[\sup_{u<R}\left|\overline{g}''_uC'_{u-}[f_j\check{g}'_{e}]
- \overline{g}''_uC''_{u-}[f_j\check{g}'_{e}]\right|]\\
&\leq&\epsilon''\mathbb{E}[\sup_{u<R}\left|C'_{u-}[f_j\check{g}'_{e}]
- C''_{u-}[f_j\check{g}'_{e}]\right|]\\

&\leq&\epsilon''\mathbb{E}[\sup_{u<R}\left|C'_{u-}[f_j\check{g}_{e}]
- C''_{u-}[f_j\check{g}_{e}]\right|]\\
&&+\epsilon''\mathbb{E}[\sup_{u<R}\left|C'_{u-}[f_j(\check{g}'_{e}-\check{g}_{e})]
\right|]
+\epsilon''\mathbb{E}[\sup_{u<R}\left|C''_{u-}[f_j(\check{g}'_{e}-\check{g}_{e})\right|]\\

&\leq&\epsilon''\|C'-C''\|_{R-}+2\epsilon''\|\check{g}'_{e}-\check{g}_{e}\|_\infty,
\dce
$$
and
$$
\dcb
&&\mathbb{E}[\sup_{u<R}\left|\overline{g}''_uC''_{u-}[f_j\check{g}'_{e}]
- \overline{g}''_uC''_{u-}[f_j\check{g}''_{e}]\right|]\\

&\leq&\epsilon''\mathbb{E}[\sup_{u<R}\left|C''_{u-}[f_j\check{g}'_{e}]
- C''_{u-}[f_j\check{g}''_{e}]\right|]\\

&\leq&\epsilon''\mathbb{E}[\sup_{u<R}\left|C''_{u-}[|\check{g}'_{e}
- \check{g}''_{e}|]\right|]\\

&\leq&\epsilon''\|\check{g}'_{e}-\check{g}''_{e}\|_\infty.
\dce
$$
In the same way, $$
\dcb
&&\mathbb{E}[\sup_{u<R}\left|\overline{g}'_uC'_{u-}[f_j]C'_{u-}[\check{g}'_{e}]
- \overline{g}''_uC'_{u-}[f_j]C'_{u-}[\check{g}'_{e}]\right|] 
\leq  \mathbb{E}[\sup_{u<R}\left|\overline{g}'_u
- \overline{g}''_u\right|]\\

&&\hspace{4cm}\leq \sum_{i=1}^k\epsilon'\|C'-C''\|_{R-}
+2\sum_{i=1}^k\epsilon'\|h_i-h'_i\|_\infty
+\sum_{i=1}^k\epsilon'\|h'_i-h''_i\|_\infty,
\\

&&\mathbb{E}[\sup_{u<R}\left|\overline{g}''_uC'_{u-}[f_j]C'_{u-}[\check{g}'_{e}]
- \overline{g}''_uC''_{u-}[f_j]C'_{u-}[\check{g}'_{e}]\right|]
\leq \epsilon''\|C'-C''\|_{R-},
\\

&&\mathbb{E}[\sup_{u<R}\left|\overline{g}''_uC''_{u-}[f_j]C'_{u-}[\check{g}'_{e}]
- \overline{g}''_uC''_{u-}[f_j]C''_{u-}[\check{g}'_{e}]\right|] 
\leq \epsilon''\|C'-C''\|_{R-},
\\

&&\mathbb{E}[\sup_{u<R}\left|\overline{g}''_uC''_{u-}[f_j]C''_{u-}[\check{g}'_{e}]
- \overline{g}''_uC''_{u-}[f_j]C''_{u-}[\check{g}''_{e}]\right|]
\leq \epsilon''\|\check{g}'_{e}-\check{g}''_{e}\|_\infty.

\dce
$$
Consequently,
$$
\dcb
&&\mathbb{E}[\sup_{u< R}|C'_{u}[f_j]-C''_{u}[f_j]|]\\
&=&\mathbb{E}[\sup_{u<R}\left|\int_0^u (C'_{s-}[f_jk'^\top_s]
- C''_{s-}[f_jk''^\top_s]) dY^{R-}_s\right|]\\

&\leq&\sum_{e=1}^da'a\left(\sum_{i=1}^k\epsilon'\|C'-C''\|_{R-}
+2\sum_{i=1}^k\epsilon'\|h_i-h'_i\|_\infty
+\sum_{i=1}^k\epsilon'\|h'_i-h''_i\|_\infty\right) \\

&&+\sum_{e=1}^da'a\left(\epsilon''\|C'-C''\|_{R-}+2\epsilon''\|\check{g}'_{e}-\check{g}_{e}\|_\infty\right) \\

&&+\sum_{e=1}^da'a\epsilon''\|\check{g}'_{e}-\check{g}''_{e}\|_\infty\\

&&+\sum_{e=1}^da'a\left(\sum_{i=1}^k\epsilon'\|C'-C''\|_{R-}
+2\sum_{i=1}^k\epsilon'\|h_i-h'_i\|_\infty
+\sum_{i=1}^k\epsilon'\|h'_i-h''_i\|_\infty\right) \\

&&+\sum_{e=1}^da'a\epsilon''\|C'-C''\|_{R-}\\

&&+\sum_{e=1}^da'a\epsilon''\|C'-C''\|_{R-}\\

&&+\sum_{e=1}^da'a\epsilon''\|\check{g}'_{e}-\check{g}''_{e}\|_\infty\\

&\leq&da'ak\epsilon'\|C'-C''\|_{R-}
+2da'a\epsilon'\sum_{i=1}^k\|h_i-h'_i\|_\infty
+da'a\epsilon'\sum_{i=1}^k\|h'_i-h''_i\|_\infty \\

&&+da'a\epsilon''\|C'-C''\|_{R-}+2a'a\epsilon''\sum_{e=1}^d\|\check{g}'_{e}-\check{g}_{e}\|_\infty\\

&&+a'a\epsilon''\sum_{e=1}^d\|\check{g}'_{e}-\check{g}''_{e}\|_\infty\\

&&+da'ak\epsilon'\|C'-C''\|_{R-}
+2da'a\epsilon'\sum_{i=1}^k\|h_i-h'_i\|_\infty
+da'a\epsilon'\sum_{i=1}^k\|h'_i-h''_i\|_\infty \\

&&+da'a\epsilon''\|C'-C''\|_{R-}\\

&&+da'a\epsilon''\|C'-C''\|_{R-}\\

&&+a'a\epsilon''\sum_{e=1}^d\|\check{g}'_{e}-\check{g}''_{e}\|_\infty\\

&=&da'a(2k\epsilon'+3\epsilon'')\|C'-C''\|_{R-}
+4da'a\epsilon'\sum_{i=1}^k\|h_i-h'_i\|_\infty
+2da'a\epsilon'\sum_{i=1}^k\|h'_i-h''_i\|_\infty\\
&&
+2a'a\epsilon''\sum_{e=1}^d\|\check{g}'_{e}-\check{g}_{e}\|_\infty
+2a'a\epsilon''\sum_{e=1}^d\|\check{g}'_{e}-\check{g}''_{e}\|_\infty.
\dce
$$
This concludes the lemma. \ok

\

\subsection{The solution in the case of continuous coefficient}

\bethe\label{existence}
Under the conditions in Section \ref{inequality} and (\ref{kDY>-1}) (in particular the functions $\mathbf {h}=(h_i:1\leq i\leq k)$ and $\check{g}_e$ being in the family $(f_j)_{j\in\mathbb{N}}$), then there exists a process of random probability measures $C_s, s\in\mathbb{R}_+$, which satisfies$$
\left\{
\dcb
dC_t&=&\int_{\mathtt{X}}k_t(\mathbf{h},C_-,x)^\top C_{t-}(dx) dY_t,\ t\in\mathbb{R}_+,\\
\\
C_0&=&\mu.
\dce
\right.
$$ 
\ethe

\textbf{Proof.}
Let $(\mathcal{D}_n)_{n\in\mathbb{N}}$ be an increasing sequence of finite $\sigma$-algebras such that $\vee_{n\in\mathbb{N}}\mathcal{D}_n=\mathcal{B}(\mathtt{X})$. Let $\mathbf{h}_n$ and $v_n$ be $\mathcal{D}_n$ measurable functions. We suppose $$
 \|\mathbf{h}-\mathbf{h}_n\|_\infty+ \|\check{g}-\check{g}\circ v_n\|_\infty<\frac{1}{2^{n}}.
$$
We note that such a setting is always possible.
We consider the solution $C^{(n)}$ of equation (\ref{cequak}) corresponding to $\mathbf{h}_n$ and $v_n$. Let $a>0$ a real number satisfying$$
(1-da'a(2k\epsilon'+3\epsilon''))>\frac{1}{2}.
$$
By Lemma \ref{ineq1}, for any stopping time $R$ such that $\|Y^{R-}_e\|_{H^\infty}<a$ for $1\leq e\leq d$,
$$
\dcb
&&(1-da'a(2k\epsilon'+3\epsilon''))\|C^{(m)}-C^{(n)}\|_{R-}\\

&\leq&4da'a\epsilon'\sum_{i=1}^k\|h_i-h_{m,i}\|_\infty
+2da'a\epsilon'\sum_{i=1}^k\|h_{m,i}-h_{n,i}\|_\infty\\
&&
+2a'a\epsilon''\sum_{e=1}^d\|\check{g}_{e}\circ v_m)-\check{g}_{e}\|_\infty
+2a'a\epsilon''\sum_{e=1}^d\|\check{g}_{e}\circ v_m)-\check{g}_{e}\circ v_n)\|_\infty\\

&\leq&4da'a\epsilon'k\frac{1}{2^m}
+2da'a\epsilon'k(\frac{1}{2^m}+\frac{1}{2^n})
+2a'a\epsilon''d\frac{1}{2^m}
+2a'a\epsilon''d(\frac{1}{2^m}+\frac{1}{2^n}).
\dce
$$
This results to $$
\lim_{n,m\rightarrow \infty}\|C^{(m)}-C^{(n)}\|_{R-}=0.
$$
From this we deduce (cf. \cite{rw}) that the sequence $C^{(n)}_s$ almost surely converges weakly to a random probability measure that we denote by $C_s$ uniformly on the interval $[0,R)$. Note then that, for any $f_j$, $$
(C^{(n)})^{R-}_t[f_j]=\mu[f_j]+\int_0^t \int_{\mathtt{X}}f_j(x)k_s(\mathbf{h}_n, C^{(n)}_-,v_n,x)^\top C^{(n)}_{s-}(dx)\ dY^{R-}_s.
$$
The expression inside the stochastic integral depends continuously on $(C^{(n)}_s: s\in[0,T))$, on $\check{g}\circ v_n$ and on $\mathbf{h}_n$, with limit$$
\int_{\mathtt{X}}f_j(x)k_s(\mathbf{h},C_-,x)^\top C_{s-}(dx).
$$ 
(Recall that $k_s(\mathbf{h}_n, C^{(n)}_-,v_n,x)$ is defined through (\ref{kC2}).) It also is uniformly bounded. Hence, the stochastic dominated convergence theorem (cf. \cite[Theorem 9.27]{Yan}) implies$$
C_t[f_j]=\mu[f_j]+\int_0^t \int_{\mathtt{X}}f_j(x)k_s(\mathbf{h}, C_-,x)^\top C_{s-}(dx)\ dY^{R-}_s,\ t<R.
$$
Consider the value at $R$ of $C^{(n)}$:
$$
C^{(n)}_{R}[f_j]=C^{(n)}_{R-}[f_j]+\Delta_RC^{(n)}[f_j]=C^{(n)}_{R-}[f_j]+\int_{\mathtt{X}}f_j(x)k_R(\mathbf{h}_n, C^{(n)}_-,v_n,x)^\top C^{(n)}_{R-}(dx)\Delta_RY. 
$$
Again, the expression on the right hand side converges to$$
C_{R-}[f_j]+\int_{\mathtt{X}}f_j(x)k_R(\mathbf{h}, C_-,x)^\top C_{R-}(dx)\Delta_RY.
$$ 
(Recall that the variable $x$ in $k_s(\mathbf{h}_n, C^{(n)}_-,v_n,x)$ is separated from the others.) Define $\Delta_RC$ to be the stochastic integral (considered as a random measure on $\mathtt{X}$) in the above expression, and set $C_t=C^{R-}_t+\Delta_RC\ind_{\{R\leq t\}}$ for $t\in[0,R]$. By Lemma \ref{raccordement}, the statement of theorem holds on the interval $[0,R]$ by the process $C$.

To continue the proof, we recall \cite[Chapter V section 3 Theorem 5]{protter} which states that, for any real number $N>0$, there exists a finite sequence of stopping times $0=R_0\leq R_1\leq R_2\leq \ldots,\leq R_q=N$ such that $\|(Y-Y^{R_i})^{R_{i+1}-}_e\|_{H^\infty}<a$ for any $0\leq i<q$. We have proved the theorem on the time interval $[0,R_1]$. Now we consider the shifted filtration $(\mathcal{F}_{R_1+t}: t\in\mathbb{R}_+)$. We consider equation (\ref{cequak}) in the shifted filtration with the shifted semimartingale $\overline{Y}_t=Y_{R_1+t}-Y_{R_1}, t\in\mathbb{R}_+$, and with the shifted initial probability measure $C_{R_1}$. The random time $\overline{R}=R_2-R_1$ is a stopping time in the shifted filtration and the norm inequality $\|\overline{Y}^{\overline{R}-}_e\|_{H^\infty}<a$ remains valid in the shifted filtration. We prove thus the theorem with a process $\overline{C}$ on the interval $[0,\overline{R}]$ with respect to the shifted setting. Define then $$
C_t=C^{R_1}_t+(\overline{C}_{(R_2-R_1)\wedge (t-R_1)^+}-C_{R_1}), \ t\in\mathbb{R}_+.
$$
We check that the above defined process $C$ satisfies the statement of the theorem on time interval $[0,R_2]$.

Continuing this process we prove the validity of the theorem on $[0,N]$. As $N>0$ is abitrary, the theorem is proved. \ok

\bethe\label{uniqueness}
Let $\mathbf{h}$ and $\mathbf{h}'$ are two vectors of continuous functions on $\mathtt{X}$ satisfying the conditions in Theorem \ref{existence}. For any stopping time $R$ satisfying $\|Y^{R-}_e\|_{H^\infty}\leq a$ for any $1\leq e\leq d$, we have
$$
\dcb
(1-da'a(2k\epsilon'+3\epsilon''))\|C'-C''\|_{R-}

&\leq&2da'a\epsilon'\sum_{i=1}^k\|h'_i-h''_i\|_\infty.
\dce
$$
Consequently, equation (\ref{cequak}) in Theorem \ref{existence} has uniqueness of solution.
\ethe

\textbf{Proof.} Repeat the argument in Lemma \ref{ineq1} and the argument based on \cite[Chapter V section 3 Theorem 5]{protter}. \ok 

\bethe\label{abscont}
Suppose that $\mathtt{X}$ is a subspace in an euclidian space. The solution to equation (\ref{cequak}) in Theorem \ref{existence} has density functions with respect to the initial probability measure $\mu$: $$
\frac{dC_t}{d\mu}(x)
=
\mathcal{E}(k(\mathbf{h},C_-,x)^\top\centerdot Y)_t,\ t\in\mathbb{R}_+, x\in\mathtt{X}.
$$
\ethe

\textbf{Proof.} Let $x_0\in\mathtt{X}$ and a real number $r>0$ be such that $\mu(B_r(x_0))>0$, where $B_r(x_0)$ be the ball at center $x_0$ of radius $r$. We write the equation
$$
\dcb
&&C_t(B_r(x_0))
=\mu(B_r(x_0))+\int_0^t \left(\int_{\mathtt{X}}\ind_{B_r(x_0)}(x)k_s(\mathbf{h},C_-,x)^\top C_{s-}(dx)\right)dY_s.
\dce
$$
or
$$
\dcb
&&\frac{C_t(B_r(x_0))}{\mu(B_r(x_0))}\\

&=&1+\int_0^t \frac{C_{s-}(B_r(x_0))}{\mu(B_r(x_0))}\left(\frac{1}{C_{s-}(B_r(x_0))}\int_{\mathtt{X}}\ind_{B_r(x_0)}(x)k_s(\mathbf{h},C_-,x)^\top C_{s-}(dx)\right)dY_s.
\dce
$$
This implies
$$
\frac{C_t(B_r(x_0))}{\mu(B_r(x_0))}=\mathcal{E}(\left(\frac{1}{C_-(B_r(x_0))}C_{-}[\ind_{B_r(x_0)}k_s(\mathbf{h},C_-)^\top]\right)\centerdot Y)_t,\ t\in\mathbb{R}_+.
$$
Because of the continity of $k_s(\mathbf{h},C_-,x)$ in $x$ and the uniform boundedness of $\frac{1}{C_-(B_r(x_0))}C_{-}[\ind_{B_r(x_0)}k_s(\mathbf{h},C_-)^\top]$, the stochastic dominated convergence theorem applies. We obtain$$
\lim_{r\downarrow 0}\frac{C_t(B_r(x_0))}{\mu(B_r(x_0))}
=
\mathcal{E}(k(\mathbf{h},C_-,x_0)^\top\centerdot Y)_t.
$$
The theorem is now the consequence of \cite[Chapter V section 4 Theorem 11]{protter} and of \cite[Theorem 4.3.4]{KP}.
\ok

\bcor\label{coroll}
The above theorems remain valid, if equation (\ref{cequak}) has a coefficient $k_t(\mathbf{h},\mathbf{c},x)$ which satisfies condition (\ref{kDY>-1}) and is a finite sum of coefficients satisfying the conditions in Theorem \ref{existence}.
\ecor

\textbf{Proof.} We need only to note that in Theorem \ref{existence} we can take a semimartingale $Y$ which has some components equal among them. \ok

\

\section{Equation (\ref{cequak}) with differentiable coefficient}

We notice that the key point in Section \ref{concoe} to establish the existence of the solution of equation (\ref{cequak}) is the inequality proved in Lemma \ref{ineq1}. It is this proof which makes us to impose the coefficient of the equation to have the particular form in subsection \ref{inequality}. In this section we will show that an inequlity like that in Lemma \ref{ineq1} can be obtained, if the coefficient is sufficiently differentiable. We establish in consequence the existence and uniqueness of equation (\ref{cequak}) with differential coefficient.

To have the differentiability we work on the space $\mathtt{X}=[0,T]^\mathfrak{l}$, where $T>0$ is a real number and $\mathfrak{l}$ is a positive integer. To make use of the differentiability we consider an initial probability measure $\mu$ which is absolutely continuous with respect to the Lebesgue measure $dx$ whose density function is denoted by $\mathfrak{m}(x)$. We consider a function $g_s(\omega,\mathbf{y},x)$ in definition (\ref{kC2}). We will say that this function is differentiable, if, in addition of its boundedness and its Lipschitzian property and condition (\ref{kDY>-1}), the function $g_s(\omega,\mathbf{y},x)$ is $\mathfrak{l}$-times differentiable with all derivatives to be Lipschitzian.

As in Section \ref{inequality}, we introduce $(f_j)_{j\in\mathbb{N}}$ a countable family of Borel functions on $\mathtt{X}$ uniformly bounded by 1, closed by multiplication. But differently to Section \ref{inequality}, we do not suppose them continuous functions. Instead, we suppose that the functions $(f_j)_{j\in\mathbb{N}}$ generate the Borel $\sigma$-algebra of $\mathtt{X}$. We suppose in particular that this family contains the functions of the form $\ind_{[0,\mathbf{b}]}$, where $\mathbf{b}=(b_1,\ldots,b_\mathfrak{l})$ with $b_i\in\mathbb{Q}_+$ and $[0,\mathbf{b}]=\prod_{i=1}^\mathfrak{l}[0,b_i]$. As in Lemma \ref{ineq1} we introduce, for any stopping time $R$, for $\mathbf{c}', \mathbf{c}''$ proceses in $\mathcal{M}(\mathtt{X})$, $$
\|\mathbf{c}'-\mathbf{c}''\|_{R-}=\sup_{j\in\mathbb{N}}\mathbb{E}[\sup_{u< R}|\mathbf{c}'_{u}[f_j]-\mathbf{c}''_{u}[f_j]|].
$$ 

\subsection{Iterated integration by parts formula}

We will need to apply the integration by parts formula repeatedly. This impose a careful organisation of the notations. This subsection is devoted to provide such a notation system.

For any integer $n>0$, let $G_n=\{1,2,\ldots,n\}$. For any $J\subset G_n$ denote $J^c_n=G_n\setminus J$ and
$
\partial_J=\prod_{j\in J}\frac{\partial}{\partial x_j}.
$

\bl\label{iterated}
Let $v,u$ be two smooth functions. We have$$
u\partial^nv
=\sum_{i=0}^n \sum_{I\subset G_n:\#I=i}(-1)^i\partial_{I^c_n}(v \partial_Iu).
$$
where $\partial^{n}=\partial_{G_{n}}$. 
\el

\textbf{Proof.}
For $n=1$, the formula is valid with $$
\dcb
u\frac{\partial v}{\partial x}=\frac{\partial}{\partial x}(uv)-v\frac{\partial u}{\partial x}.
\dce
$$
And, applying the first formula, we obtain a second formula:
$$
\dcb
u\frac{\partial}{\partial y}\frac{\partial v}{\partial x}
&=&
\frac{\partial}{\partial y}(u\frac{\partial v}{\partial x})-\frac{\partial u}{\partial y}\frac{\partial v}{\partial x}\\

&=&
\frac{\partial}{\partial y}(\frac{\partial}{\partial x}(uv)-\frac{\partial u}{\partial x}v)-\left(\frac{\partial}{\partial x}(\frac{\partial u}{\partial y}v)-\frac{\partial^2 u}{\partial y\partial x}v\right)\\

&=&
\frac{\partial}{\partial y\partial x}(uv)-\frac{\partial}{\partial y}(\frac{\partial u}{\partial x}v)-\frac{\partial}{\partial x}(\frac{\partial u}{\partial y}v)+\frac{\partial^2 u}{\partial y\partial x}v.
\dce
$$
Rearranging this formula we obtain
$$
\dcb
u\frac{\partial^2 v}{\partial x_2\partial x_1}

&=&
\frac{\partial}{\partial x_2\partial x_1}(uv)-\frac{\partial}{\partial x_2}(\frac{\partial u}{\partial x_1}v)-\frac{\partial}{\partial x_1}(\frac{\partial u}{\partial x_2}v)+\frac{\partial^2 u}{\partial x_2\partial x_1}v\\

&=&
\sum_{i=0}^2 \sum_{I\subset G_2:\#I=i}(-1)^i\partial_{I^c_2}(\partial_Iu \ v).
\dce
$$
Suppose that the above formula holds till the level $n-1$ ($n\in\mathbb{N}^*$):$$
u\partial^{n-1}v
=
\sum_{i=0}^{n-1} \sum_{I\subset G_{n-1}:\#I=i}(-1)^i\partial_{I^c_{n-1}}(\partial_Iu \ v).
$$
For the level $n$, we have
$$
\dcb
u\partial^nv
&=&u\partial^{n-1}\frac{\partial v}{\partial x_n}\\

&=&\sum_{i=0}^{n-1} \sum_{I\subset G_{n-1}:\#I=i}(-1)^i\partial_{I^c_{n-1}}(\partial_Iu \ \frac{\partial v}{\partial x_n})\\

&=&\sum_{i=0}^{n-1} \sum_{I\subset G_{n-1}:\#I=i}(-1)^i\partial_{I^c_{n-1}}(\frac{\partial}{\partial x_n}(\partial_Iu\ v)- \frac{\partial }{\partial x_n}\partial_Iu\ v)\\

&=&\sum_{i=0}^{n-1} \sum_{I\subset G_{n-1}:\#I=i}(-1)^i\partial_{I^c_{n-1}}\frac{\partial}{\partial x_n}(\partial_Iu\ v)
-
\sum_{i=0}^{n-1} \sum_{I\subset G_{n-1}:\#I=i}(-1)^i\partial_{I^c_{n-1}}(\frac{\partial }{\partial x_n}\partial_Iu\ v)\\

&=&\sum_{i=0}^{n-1} \sum_{I\subset G_{n}:\#I=i, n\notin I}(-1)^i\partial_{I^c_{n}}(\partial_Iu\ v)
-
\sum_{i=0}^{n-1} \sum_{I\subset G_{n}:\#I=i+1,n\in I}(-1)^i\partial_{I^c_{n}}(\partial_Iu\ v)\\

&=&\sum_{i=0}^{n-1} \sum_{I\subset G_{n}:\#I=i, n\notin I}(-1)^i\partial_{I^c_{n}}(\partial_Iu\ v)
+
\sum_{i=1}^{n} \sum_{I\subset G_{n}:\#I=i,n\in I}(-1)^i\partial_{I^c_{n}}(\partial_Iu\ v)\\

&=&
\sum_{i=0}^n \sum_{I\subset G_n:\#I=i}(-1)^i\partial_{I^c_n}(\partial_Iu \ v).
\dce
$$
\ok

\bl\label{integparts}
For real numbers $0\leq a_{j,0}<a_{j,1}, j\in G_n$, let $D=\prod_{j\in G_n}(a_{j,0},a_{j,1}]$. For a bounded Borel function $f$, let$$
F(x)=
\int_0^{x_1}\ldots\int_0^{x_{n}}f(z)dz.
$$
Then, ,
$$
\int_{D}f(x) dx
=
\sum_{\epsilon_j\in\{0,1\}: j\in G_n}(-1)^{\sum_{j\in G_n} (1-\epsilon_j)}F(a_{j,\epsilon_j}:j\in G_n).
$$
\el

\textbf{Proof.}
Notice that the above expression is the so-called $\mathfrak{l}$-volume (cf. \cite[Definition 2.10.1]{nelsen}). Let us prove the formula of the lemma when $f$ is smooth. The general case will then be the consequence of the monotone class theorem. For $n=1$,
$$
\int_{a_{1,0}}^{a_{1,1}}\partial F(x_1) dx_1
=
F({a_{1,1}})-F({a_{1,0}})
=
\sum_{\epsilon_j\in\{0,1\}: j\in G_1}(-1)^{\sum_{j\in G_1} (1-\epsilon_j)}F(a_{j,\epsilon_j}:j\in G_1).
$$
Suppose that the formula remains valid till level $n-1$. At the level $n$, we have
$$
\dcb
&&\int_{D_n}\partial^n F(x) dx\\
&=&\int_{a_{n,0}}^{a_{n,1}}\int_{D_{n-1}}\partial_{x_n}\partial^{n-1} F(x_{G_{n-1}},x_n) dx_{G_{n-1}}dx_{n}\\

&=&\int_{a_{n,0}}^{a_{n,1}}\left(\sum_{\epsilon_j\in\{0,1\}: j\in G_{n-1}}(-1)^{\sum_{j\in G_{n-1}} (1-\epsilon_j)}\partial_{x_n}F((a_{j,\epsilon_j}:j\in G_{n-1}),x_n)\right)dx_{n}\\

&=&\sum_{\epsilon_j\in\{0,1\}: j\in G_{n-1}}(-1)^{\sum_{j\in G_{n-1}} (1-\epsilon_j)}F((a_{j,\epsilon_j}:j\in G_{n-1}),a_{n,1})\\
&&-\sum_{\epsilon_j\in\{0,1\}: j\in G_{n-1}}(-1)^{\sum_{j\in G_{n-1}} (1-\epsilon_j)}F((a_{j,\epsilon_j}:j\in G_{n-1}),a_{n,0})\\

&=&\sum_{\epsilon_j\in\{0,1\}: j\in G_{n},\epsilon_n=1}(-1)^{\sum_{j\in G_{n}} (1-\epsilon_j)}F(a_{j,\epsilon_j}:j\in G_{n})\\
&&+\sum_{\epsilon_j\in\{0,1\}: j\in G_{n},\epsilon_n=0}(-1)^{\sum_{j\in G_{n}} (1-\epsilon_j)}F(a_{j,\epsilon_j}:j\in G_n)\\

&=&\sum_{\epsilon_j\in\{0,1\}: j\in G_{n}}(-1)^{\sum_{j\in G_{n}} (1-\epsilon_j)}F(a_{j,\epsilon_j}:j\in G_n).
\dce
$$
\ok

\bl\label{fFu}
For the same functions $f, F$ as above computed in the case $n=\mathfrak{l}$, for $\mathfrak{l}$-times continuously differentiable function $u(x)$, $$
\dcb
&&
\int_{\mathtt{X}}f(x) u(x)dx\\

&=&\sum_{i=0}^\mathfrak{l} \sum_{I\subset G_\mathfrak{l}:\#I=i}(-1)^i\int_{[0,T]^I}F((T:j\in I_\mathfrak{l}^c),x_I)\partial_Iu((T:j\in I_\mathfrak{l}^c),x_I) dx_I.
\dce
$$
\el

\textbf{Proof.}
We need only to prove the lemma for smooth function $f$. We apply Lemma \ref{iterated} and \ref{iterated} to write$$
\dcb
&&
\int_{\mathtt{X}}\partial^{\mathfrak{l}}F(x) u(x)dx\\

&=&\int_{\mathtt{X}}\sum_{i=0}^\mathfrak{l} \sum_{I\subset G_\mathfrak{l}:\#I=i}(-1)^i\partial_{I^c_n}(F\partial_Iu)(x) dx\\

&=&\sum_{i=0}^\mathfrak{l} \sum_{I\subset G_\mathfrak{l}:\#I=i}(-1)^i\sum_{\epsilon_j\in\{0,1\}: j\in I_\mathfrak{l}^c}(-1)^{\sum_{j\in I_\mathfrak{l}^c} (1-\epsilon_j)}\\
&&\int_{[0,T]^I}F((\epsilon_jT:j\in I_\mathfrak{l}^c),x_I)\partial_Iu((\epsilon_jT:j\in I_\mathfrak{l}^c),x_I) dx_I\\

&=&\sum_{i=0}^\mathfrak{l} \sum_{I\subset G_\mathfrak{l}:\#I=i}(-1)^i\int_{[0,T]^I}F((T:j\in I_\mathfrak{l}^c),x_I)\partial_Iu((T:j\in I_\mathfrak{l}^c),x_I) dx_I.
\dce
$$
\ok

\subsection{The inequality}

We now consider two differentiable functions $g'_s(\omega,\mathbf{y},x)$ and $g''_s(\omega,\mathbf{y},x)$ and two vectors $\mathbf{h}'$ and $\mathbf{h}''$ of bounded Borel functions on $\mathtt{X}$. Consider the corresponding equation (\ref{cequak}) with initial probability measure $\mu(dx)=\mathfrak{m}(x)dx$ and suppose that they have solutions respectively $C',C''$, which are absolutely continuous with respect to $\mu$.  
Denote $\frac{dC'_{u-}}{d\mu}(x)=\mathsf{p}'_{u-}(x)$ and $\frac{dC''_{u-}}{d\mu}(x)=\mathsf{p}''_{u-}(x)$ and$$
\dcb
g'_s(x)=g'(\omega,C'[\mathbf{h}'],x), \ \mbox{ and similarly $g''_s(x)$},\\

k'_s=k_s(\omega,\mathbf{h}',C'_-,x),\ \mbox{ and similarly $k''_s$}.
\dce
$$
We denote uniformly all the absolute bound of the derivatives $\partial_Jg'_s$ and $\partial_Jg''_s$, $J\subset G_\mathfrak{l}$, by $\epsilon''$.

\bl\label{ineq2}
Let $a'$ be a constant in \cite[Chapter V section 2 Theorem 2]{protter}. Suppose that the stopping time $R$ satisfies $\|Y^{R-}_e\|_{H^\infty}\leq a$, for some real $a>0$ and for any $1\leq e\leq d$. Then,
$$
\dcb
&&(1-da'a(2\epsilon''(T+1)^\mathfrak{l}+\epsilon''))\|C'-C''\|_{R-}

\leq
2\sum_{e=1}^da'a\mathbb{E}[\sup_{u<R}\|g'_{e,u}-g''_{e,u}\|_\infty].
\dce
$$
\el

\textbf{Proof.} We know that to prove this inequality basically is to make a good estimation on the right hand side term of the followng computations: $$
\dcb
&&\mathbb{E}[\sup_{u<R}\left|C'_{u-}[f_jk'_{e,u}]
- C''_{u-}[f_jk''_{e,u}]\right|]\\

&\leq&
\mathbb{E}[\sup_{u<R}\left|C'_{u-}[f_jg'_{e,u}]
- C''_{u-}[f_jg''_{e,u}]\right|]
+\mathbb{E}[\sup_{u<R}\left|C'_{u-}[f_j]C'_{u-}[g'_{e,u}]
- C''_{u-}[f_j]C''_{u-}[g''_{e,u}]\right|]\\

&\leq&
\mathbb{E}[\sup_{u<R}\left|C'_{u-}[f_jg'_{e,u}]
- C''_{u-}[f_jg'_{e,u}]\right|]
+\mathbb{E}[\sup_{u<R}\left|C''_{u-}[f_jg'_{e,u}]
- C''_{u-}[f_jg''_{e,u}]\right|]\\
&&+\mathbb{E}[\sup_{u<R}\left|C'_{u-}[f_j]C'_{u-}[g'_{e,u}]
- C''_{u-}[f_j]C'_{u-}[g'_{e,u}]\right|]
+\mathbb{E}[\sup_{u<R}\left|C''_{u-}[f_j]C'_{u-}[g'_{e,u}]
- C''_{u-}[f_j]C''_{u-}[g'_{e,u}]\right|]\\
&&+\mathbb{E}[\sup_{u<R}\left|C''_{u-}[f_j]C''_{u-}[g'_{e,u}]
- C''_{u-}[f_j]C''_{u-}[g''_{e,u}]\right|]\\

&\leq&
\mathbb{E}[\sup_{u<R}\left|C'_{u-}[f_jg'_{e,u}]
- C''_{u-}[f_jg'_{e,u}]\right|]
+\mathbb{E}[\sup_{u<R}\|g'_{e,u}-g''_{e,u}\|_\infty]\\
&&+\epsilon''\mathbb{E}[\sup_{u<R}\left|C'_{u-}[f_j]
- C''_{u-}[f_j]\right|]
+\mathbb{E}[\sup_{u<R}\left|C'_{u-}[g'_{e,u}]
- C''_{u-}[g'_{e,u}]\right|]\\
&&+\mathbb{E}[\sup_{u<R}\|g'_{e,u}-g''_{e,u}\|_\infty].

\dce
$$
We introduce $$
\dcb
v'_u(x)
=
\int_0^{x_1}\ldots\int_0^{x_{\mathfrak{l}}}f_j(x)\mathsf{p}'_{u-}(x) \mathfrak{m}(x) dx=C'_{u-}[f_j\ind_{[0,x]}],
\dce
$$
and in the same way $v''_t(x)$. Then, by Lemma \ref{fFu},
$$
\dcb
&&C'_{u-}[f_jg'_{e,u}]
- C''_{u-}[f_jg'_{e,u}]\\

&=&
\int_{\mathtt{X}}\left(f_j(x)\mathsf{p}'_{u-}(x)\mathfrak{m}(x)
-f_j(x)\mathsf{p}''_{u-}(x)\mathfrak{m}(x)\right) g'_{e,u}(x)dx\\

&=&\sum_{i=0}^\mathfrak{l} \sum_{I\subset G_\mathfrak{l}:\#I=i}(-1)^i\int_{[0,T]^I}(v'_u-v''_u)((T:j\in I_\mathfrak{l}^c),x_I)\partial_Ig_u((T:j\in I_\mathfrak{l}^c),x_I) dx_I\\

&=&\sum_{i=0}^\mathfrak{l} \sum_{I\subset G_\mathfrak{l}:\#I=i}(-1)^i\int_{[0,T]^I} \left(C'_{u-}[f_j\ind_{[,\mathbf{b}_I(x_I)]}]-C''_{u-}[f_j\ind_{[,\mathbf{b}_I(x_I)]}]\right)\partial_Ig_u((T:j\in I_\mathfrak{l}^c),x_I) dx_I,
\dce
$$
where $\mathbf{b}_I(x_I)$ is the vector whose component at $i\in I$ is given by $x_i$ and that at $i\notin I$ given by $T$. This identity implies the estimation
$$
\dcb
&&\mathbb{E}[\sup_{u<R}\left|C'_{u-}[f_jg'_{e,u}]
- C''_{u-}[f_jg'_{e,u}]\right|]\\

&\leq&\sum_{i=0}^\mathfrak{l} \sum_{I\subset G_\mathfrak{l}:\#I=i}\epsilon''\int_{[0,T]^I}\mathbb{E}[\sup_{u<R}\left|C'_{u-}[f_j\ind_{[,\mathbf{b}_I(x_I)]}]-C''_{u-}[f_j\ind_{[,\mathbf{b}_I(x_I)]}]\right|]dx_I\\

&\leq&\sum_{i=0}^\mathfrak{l} \sum_{I\subset G_\mathfrak{l}:\#I=i}\epsilon''\int_{[0,T]^I}\|C'-C''\|_{R-}dx_I\\

&\leq&\epsilon''(T+1)^\mathfrak{l}\|C'-C''\|_{R-}.
\dce
$$
Applying the above estimation with $f_j\equiv 1$, we also obtain
$$
\dcb
&&\mathbb{E}[\sup_{u<R}\left|C'_{u-}[g'_{e,u}]
- C''_{u-}[g'_{e,u}]\right|]

\leq \epsilon''(T+1)^\mathfrak{l}\|C'-C''\|_{R-}.
\dce
$$
Together with the obvious inequality$$
\mathbb{E}[\sup_{u<R}\left|C'_{u-}[f_j]
- C''_{u-}[f_j]\right|]
\leq \|C'-C''\|_{R-},
$$
we conclude$$
\mathbb{E}[\sup_{u<R}\left|C'_{u-}[f_jk'_{e,u}]
- C''_{u-}[f_jk''_{e,u}]\right|]
\leq
(2\epsilon''(T+1)^\mathfrak{l}+\epsilon'')\|C'-C''\|_{R-}
+2\mathbb{E}[\sup_{u<R}\|g'_{e,u}-g''_{e,u}\|_\infty].
$$
Now we can write
$$
\dcb
&&\mathbb{E}[\sup_{u< R}|C'_{u}[f_j]-C''_{u}[f_j]|]\\

&=&\mathbb{E}[\sup_{u< R}|C'^{R-}_{u}[f_i]-C''^{R-}_{u}[f_i]|]\\

&=&\mathbb{E}[\sup_{u<R}\left|\int_0^u (C'_{s-}[f_jk'^\top_s]
- C''_{s-}[f_jk''^\top_s]) dY^{R-}_s\right|]\\

&\leq&\sum_{e=1}^d\mathbb{E}[\sup_{u<R}\left|\int_0^u (C'_{s-}[f_jk'_{e,s}]
- C''_{s-}[f_jk''_{e,s}]) dY^{R-}_{e,s}\right|]\\

&\leq&\sum_{e=1}^da'\mathbb{E}[\sup_{u<R}\left|C'_{u-}[f_jk'_{e,u}]
- C''_{u-}[f_jk''_{e,u}]\right|]\|Y^{R-}_{e}\|_{H^\infty}\\

&\leq&\sum_{e=1}^da'a\mathbb{E}[\sup_{u<R}\left|C'_{u-}[f_jk'_{e,u}]
- C''_{u-}[f_jk''_{e,u}]\right|]\\

&\leq&\sum_{e=1}^da'a\left((2\epsilon''(T+1)^\mathfrak{l}+\epsilon'')\|C'-C''\|_{R-}
+2\mathbb{E}[\sup_{u<R}\|g'_{e,u}-g''_{e,u}\|_\infty]\right)\\

&\leq&da'a(2\epsilon''(T+1)^\mathfrak{l}+\epsilon'')\|C'-C''\|_{R-}
+2a'a\sum_{e=1}^d\mathbb{E}[\sup_{u<R}\|g'_{e,u}-g''_{e,u}\|_\infty].

\dce
$$
\ok

\subsection{The solution in the case of differentiable coefficient}

Let $g_s(\omega,\mathbf{y},x)$ be a differential function and $\mathbf{h}$ be a vector of continuous functions. Let $k_s(\mathbf{h},\mathbf{c},x)$ be the associated coefficient in (\ref{kC2}). Suppose that there exists a sequence of differential functions $g_{n,s}(\omega,\mathbf{y},x), n\in\mathbb{N}$, which in addition satisfy the conditions Corollary \ref{coroll}, such that
$$
\sum_{I\subset G_\mathfrak{l}}\sum_{e=1}^d\sup_{u<R}\sup_{\omega,\mathbf{y},x}|\partial_Ig_{e,u}(\omega,\mathbf{y},x)-\partial_Ig_{n,e,u}(\omega,\mathbf{y},x)|_\infty\leq \frac{1}{2^n}.
$$

\bethe\label{existence2}
Under the above condition, there exists a unique process of random probability measures $C_s, s\in\mathbb{R}_+$, which satisfies$$
\left\{
\dcb
dC_t&=&\int_{\mathtt{X}}k_t(\mathbf{h},C_-,x)^\top C_{t-}(dx) dY_t,\ t\in\mathbb{R}_+,\\
\\
C_0&=&\mu.
\dce
\right.
$$ 
The solutions $C_t$ have density functions with respect to the initial probability measure $\mu$: $$
\frac{dC_t}{d\mu}(x)
=
\mathcal{E}(k(\mathbf{h},C_-,x)^\top\centerdot Y)_t,\ t\in\mathbb{R}_+, x\in\mathtt{X}.
$$
\ethe

\textbf{Proof.} Note that under the above condition, $\partial_Ig_{n,e,u}$ is uniformly bounded. The inequality in Lemma \ref{ineq2} can be used with a same constant $\epsilon''$. We now only need to repeat the argument in the proof of Theorem \ref{existence} to establish the existence of equation (\ref{cequak}) in this theorem. In particular, considering equations (\ref{cequak}) corresponding to the functions $g_n$, with the same vector $\mathbf{h}$ and the same initial probability measure $\mathfrak{m}(x)dx$, the solutions $C^{(n)}, n\in\mathbb{N}$ of these equations converge to the solution of the equation in this theorem, at least on an non empty interval $[0,R]$. \ok

\

\end{document}